\documentclass[a4paper,11pt,twoside]{amsart}
\usepackage[latin1]{inputenc}
\usepackage[T1]{fontenc}

\usepackage{a4wide}

\usepackage{ifpdf}
\ifpdf\usepackage[pdftex]{hyperref}
\else\usepackage[hypertex]{hyperref}\fi

\usepackage{amsmath,amssymb,amsthm}

\usepackage{lmodern}

\theoremstyle{plain}
\newtheorem{thm}{Theorem}[section]
\newtheorem{prop}[thm]{Proposition}
\newtheorem{lemma}[thm]{Lemma}
\newtheorem{cor}[thm]{Corollary}

\theoremstyle{definition}
\newtheorem{defn}[thm]{Definition}

\theoremstyle{remark}
\newtheorem{rem}[thm]{Remark}

\theoremstyle{remark}

\DeclareMathOperator{\Sp}{Sp}

\def\R{\mathbb R}

\def\N{\mathbb N}
\def\Z{\mathbb Z}

\begin{document}

\title[Uncertainty principle for POVM]{An entropic uncertainty principle\\
  for positive operator valued measures}

\author{Michel Rumin}
\address{Laboratoire de Math\'ematiques d'Orsay\\
  UMR 8628 \\ CNRS et Universit\'e Paris-Sud\\
  91405 Orsay\\ France}

\email{michel.rumin@math.u-psud.fr}

\date{\today}

\begin{abstract}
  Extending a recent result by Frank and Lieb, we show an entropic
  uncertainty principle for mixed states in a Hilbert space, relatively to
  pairs of positive operator valued measures that are independent in some
  sense. This yields spatial-spectral uncertainty principles and
  log-Sobolev inequalities for invariant operators on homogeneous spaces,
  which are sharp in the compact case.
\end{abstract}

\keywords{entropy, POVM, uncertainty principle, homogeneous spaces,
  log-Sobolev}

\subjclass[2000]{58J42, 58J50, 47B06, 43A85, 94A17.}


\maketitle


\section{Introduction and main result}
\label{sec:introduction}

A mixed state (or density matrix) in a Hilbert space $H$ is a positive
operator $\rho$ on $H$ with trace $\tau(\rho)= 1$. From the quantum
view-point, one can measure $\rho$ through positive operator valued
measures (POVM) on $H$. We briefly recall this.

A POVM is a countably additive map $P$ from a $\sigma$-algebra on a space
$X$ into the positive operators on $H$ and such that $P(X) = 1_H $, see
\cite{Nielsen-Chuang,Wikipedia-POVM}. Important examples are given by the
following cases:
\begin{itemize}
\item When $H$ is isometric to a function space $L^2(X,\mu)$, one can set
  $P(\Omega)$ to be the multiplication by the characteristic function
  $\chi_\Omega$ for measurable sets $\Omega \subset X$. This is actually a
  \emph{projection} valued measure.
\item Other projection valued measures come from the spectral resolution of
  a self-adjoint operator $A$ on $H$. Here $X= \R$ and $P(I)= \Pi_A(I)$ is
  the spectral projector of $A$ associated to $I \subset\R$.
\item Partition of unity in $H$, i.e.~ series of positive operators $P_i$
  such that $\sum_i P_i = 1_H$.
\end{itemize}

A state $\rho$ and a POVM induce a probability measure $\nu_P^\rho$ on $X$
by
$$
\nu_P^\rho(\Omega) = \tau(\rho P(\Omega)) = \tau(\rho^{1/2} P(\Omega)
\rho^{1/2})\,,
$$
that encodes the result of the measurement of $\rho$ through $P$.

The statement we want to discuss deals with the measures $\nu^\rho_P$ and
$\nu^\rho_Q$ associated to a \emph{pair} of POVM on $H$, say $P$ from $X$
and $Q$ from $Y$. This pair induces a state-independent measure $\mu_{PQ}$
on $X \times Y$, called a 'Liouville measure' from now on, and defined by
\begin{displaymath}
  \mu_{PQ}(\Omega_1
  \times \Omega_2) = \tau(P(\Omega_1)Q(\Omega_2)) = \tau(P(\Omega_1)^{1/2}
  Q(\Omega_2) P(\Omega_1)^{1/2}) \,.
\end{displaymath}


We shall see that when the Liouville measure $\mu_{PQ}$ is bounded by a
$\sigma$-finite product measure $\mu_P\otimes \mu_Q$, then $\nu_P^\rho \ll
\mu_P$ and $\nu_Q^\rho \ll \mu_Q $. Let then
\begin{displaymath}
  \varrho_P = \frac{d \nu_P^\rho}{d\mu_P} \quad \mathrm{and} \quad 
  \varrho_Q = \frac{d \nu_Q^\rho}{d\mu_Q}
\end{displaymath}
denote the corresponding density functions. Our main result is the
following.

\begin{thm}\label{thm:1}
  Using the notations above, suppose that $\mu_{PQ} \leq \mu_P \otimes
  \mu_Q$ for some $\sigma$-finite product measure on $X \times Y$.
  Then it holds
  \begin{equation}
    \label{eq:1} 
    -\int_X \varrho_P \ln\varrho_P d \mu_P -\int_Y \varrho_Q \ln \varrho_Q
    d \mu_Q \geq S(\rho) = -\tau(\rho \ln \rho)\,,
  \end{equation}
  provided the left side is well defined, i.e.~we don't face
  $\infty-\infty$ in the integrals and their sum.
\end{thm}

Here $S(\rho)$ is the intrinsic von Neumann entropy of the state,
independent on any choice of POVM, while the integrals are minus the
relative entropy from $\nu_P^\rho$ to $\mu_P$ and $\nu_Q^\rho$ to
$\mu_Q$. Hence \eqref{eq:1} is an entropic uncertainty principle claiming
that the $P$ and $Q$-measurements of $\rho$ can't be \emph{both}
concentrated in small sets for $\mu_P$ and $\mu_Q$, under some hypothesis
on the Liouville measure. From the probabilistic view-point, the assumption
on $\mu_{PQ}$ deals with the independence of the two POVM. The best case
occurs when $\mu_{PQ}$ is a product measure itself: and we shall say that
$P$ are $Q$ are \emph{independent} there. We give some examples.

\section{First examples}
\label{sec:first-examples}


\subsection{Pair of bases}
\label{sec:pair-bases}

Let $(e_i)$ and $(f_i)$ be two orthonormal bases of a Hilbert space
$H$. This gives two projector valued measures from (finite or not) sets
$X=Y \subset \Z$ with $P(i) = \Pi_{e_i}$ and $Q(i) = \Pi_{f_i}$.  The
measurements of $\rho$ are
$$
\nu_P^\rho (i) = \langle \rho e_i, e_i \rangle = p_i \ , \quad \nu_Q^\rho
(j) = \langle \rho f_j , f_j \rangle = q_j\,
$$ 
and the Liouville measure is
\begin{displaymath}
  \mu_{PQ}(i,j) = \tau (P(i) Q(j)) = |\langle e_i, f_j\rangle|^2\,.
\end{displaymath}
Clearly, it holds that $\mu_{PQ} \leq \mu_P \otimes \mu_Q$ with
\begin{displaymath}
  \mu_P(i) = \sup_{j}|\langle e_i, f_j\rangle| \quad \mathrm{and}
  \quad \mu_Q(j) = \sup_{i}|\langle e_i, f_j\rangle | .
\end{displaymath}
Theorem~\ref{thm:1} yields
\begin{equation}
  \label{eq:2}
  \begin{gathered}
    -\sum_i p_i \ln p_i  - \sum_j  q_j \ln q_j \geq S(\rho) + K \,,\\
    \mathrm{with} \quad K = - \sum_i p_i \ln \mu_P(i) - \sum_j q_j \ln
    \mu_Q(j) \geq 0 \,.
  \end{gathered}
\end{equation}
This improves a result by Frank and Lieb in \cite{Frank-Lieb} where $K$ is
replaced by
$$
K'= -2 \ln (\sup_{i,j} |\langle e_i, f_j \rangle|) \leq K .
$$

The better estimation is obtained when $P$ and $Q$ are independent, i.e.~
when $\mu_{PQ}$ is a product measure. Taking partial $i$ and $j$ sums, one
sees that this only happens in finite dimension $d$ and when
\begin{displaymath}
  |\langle e_i, f_j \rangle |^2 = d^{-1} \,,
\end{displaymath}
which is a definition of \emph{mutually unbiased bases}, used in quantum
information theory \cite{Nielsen-Chuang,Wikipedia-MUB}.

\smallskip We note that such results generalize to pairs of partitions of
unity $P_i$ and $Q_j$ obtained in the following way. Given two orthonormal
bases $(e_i)$ and $f_j$ of Hilbert space $H'$ containing $H$, one gets two
POVM on $H$ by setting
$$
P_i = \Pi_H \Pi_{e_i} \Pi_H \quad \mathrm{and}\quad Q_j= \Pi_H \Pi_{f_j}
\Pi_H \,.
$$ 
One finds that the same inequality as \eqref{eq:2} holds replacing $e_i$
and $f_j$ by $\Pi_H e_i$ and $\Pi_H f_j$ in formulas.

\subsection{Fourier transform and uncertainty}
\label{sec:fourier-uncertainty}

Consider now $H= L^2(\R^n, dx)$ with its natural projection valued measure
from $X=\R^n_x$ defined by $P(\Omega) = \chi_\Omega \times$.  Fourier
transform on $H$, with convention that
$$
\mathcal{F} (f) (\xi) = \int_{\R^n} f(x) e^{-2i\pi\langle x,\xi\rangle} dx
$$
is an isometry onto $L^2(\R^n_\xi, d \xi)$, and provides a second
projection valued measure $Q$ from $Y= \R^n_\xi$ into $H$ by setting
$Q(\Omega) = \mathcal{F}^{-1}\chi_\Omega \mathcal{F} $.

These two maps are independent in our sense. Indeed, our 'Liouville
measure' is here
\begin{align*}
  \mu_{PQ}(\Omega_1 \times \Omega_2) & = \tau( \chi_{\Omega_1}
  \mathcal{F}^* \chi_{\Omega_2} \mathcal{F}) = \|
  \chi_{\Omega_2}\mathcal{F} \chi_{\Omega_1}\|_{HS}^2 \
  \mathrm{(Hilbert\ Schmidt\ norm)} \\
  & = \int_{\Omega_2 \times \Omega_1} |e^{-2i \pi\langle x,\xi\rangle}|^2
  d  \xi dx \quad \mathrm{(using\ kernels)}\\
  &= \lambda (\Omega_1) \lambda(\Omega_2)
\end{align*}
i.e.~ the (genuine) Liouville measure of the phase space $\R^n_x \times
\R^n_\xi$.

We make explicit Theorem~\ref{thm:1}. A state may be written $\rho = \sum_i
p_i \Pi_{f_i}$ for orthonormal functions $f_i$. Then its measurement in $P$
reads
\begin{displaymath}
  \nu_P^\rho(\Omega) = \tau(\rho \chi_\Omega) = \int_\Omega
  \sum_i p_i |f_i(x)|^2 dx \,,
\end{displaymath} 
while from $Q$
\begin{align*}
  \nu_Q^\rho(\Omega) & = \tau(\rho \mathcal{F}^{-1} \chi_\Omega
  \mathcal{F}) = \tau(\mathcal{F} \rho \mathcal{F}^{-1} \chi_\Omega)
  \\
  & = \int_\Omega \sum_i p_i |\mathcal{F}(f_i)(\xi)|^2 d\xi \,,
\end{align*}
since $\mathcal{F} \rho \mathcal{F}^{-1} = \sum_i p_i
\Pi_{\mathcal{F}(f_i)}$. Finally, Theorem~\ref{thm:1} states that
\begin{equation}
  \label{eq:3}
  - \int_{\R^n} \varrho(x) \ln \varrho(x) d x - \int_{\R^n}
  \widehat\varrho(\xi) \ln  \widehat\varrho (\xi) d \xi  \geq S(\rho) = -
  \sum_i p_i \ln p_i
\end{equation}
with
$$
\varrho(x) = \frac{d \nu_P^\rho}{dx} = \sum_i p_i |f_i(x)|^2 \quad
\mathrm{and}\quad \widehat\varrho (\xi) = \frac{d \nu_Q^\rho}{d\xi} =
\sum_i p_i |\mathcal{F}(f_i)(\xi)|^2.
$$  
This inequality has been first proved by Frank and Lieb in
\cite{Frank-Lieb}, and is known to be sharp for states $\rho_t = e^{-t
  (\Delta + \|x\|^2)}/ \tau(e^{-t( \Delta + \|x\|^2)})$ when $t \searrow
0$.


\subsection{Frank and Lieb result}
\label{sec:frank-lieb-result}

The previous example is actually a particular case of a more general
theorem proved in \cite{Frank-Lieb}.

Suppose that two $\sigma$-finite function spaces $L^2(X,\mu)$ and $L^2(Y,
\nu)$ are isometric through
$$
\mathcal{U}: L^2(X,\mu) \rightarrow L^2(Y, \nu) \,.
$$ 
Suppose moreover that $\mathcal{U}$ is bounded from $L^1(X,\mu)$ to
$L^\infty(Y,\nu)$. Let $\rho = \sum_i p_i \Pi_{f_i}$ be a (unit trace)
state on $L^2(X,\mu)$ and consider $\widehat \rho= \mathcal{U} \rho
\mathcal{U}^{-1} = \sum_i p_i \Pi_{\mathcal{U}f_i}$ on $L^2(Y, \nu)$. Let
\begin{displaymath}
  \varrho(x) =
  \sum_i p_i |f_i(x)|^2 \quad \mathrm{and} \quad \widehat\varrho (\xi)
  = \sum_i p_i  |\mathcal{U}(f_i)(\xi)|^2
\end{displaymath}
be their density as above.

\begin{thm}\cite[Thm~2.2]{Frank-Lieb}
  \label{thm:2}
  Suppose moreover that
  \begin{displaymath}
    \int_X \varrho(x) \ln^+ \varrho(x) d \mu(x) < + \infty \quad
    \mathrm{and} \quad   \int_Y \widehat\varrho(y) \ln^+
    \widehat\varrho(y) d\nu (y) < + \infty
  \end{displaymath}
  Then
  \begin{equation}
    \label{eq:4}
    - \int_X \varrho(x) \ln \varrho(x) d \mu(x) - \int_Y
    \widehat\varrho(y) \ln \widehat \varrho(y) d\nu (y) \geq S(\rho)
    - 2 \ln \|\mathcal{U}\|_{L^1 \rightarrow L^\infty}\,.  
  \end{equation}
\end{thm}

Theorem~\ref{thm:1} implies it. Indeed, as in the previous discussion on
Fourier transform, one defines two projections valued measures on $H=
L^2(X,\mu)$: the standard one $P(\Omega) = \chi_\Omega \times$ from $X$,
and $Q(\Omega)= \mathcal{U} \chi_\Omega \mathcal{U}^{-1}$ from $Y$. Then
arguing as previously, one has
\begin{align*}
  \mu_{PQ}(\Omega_1 \times \Omega_2) & = \tau( \chi_{\Omega_1} \mathcal{U}
  \chi_{\Omega_2} \mathcal{U}^*) = \|
  \chi_{\Omega_2}\mathcal{U}^* \chi_{\Omega_1}\|_{HS}^2 \\
  & = \int_{\Omega_1 \times \Omega_2} |\mathcal{U}(y,x)|^2 d\mu(x)
  d\nu(y) \\
  & \leq \mu(\Omega_1) \nu(\Omega_2) \|\mathcal{U}\|^2_{L^1 \rightarrow
    L^\infty} \,.
\end{align*}
This shows that Theorem~\ref{thm:1} applies and yields \eqref{eq:4}.

\smallskip

Comparing to Theorem~\ref{thm:2}, Theorem~\ref{thm:1} adds some flexibility
in the choice of the mesures on $X$ and $Y$, and applies to general
positive operator valued measures, as illustrated in
\S\ref{sec:pair-bases}. Here the mesures on $X$ and $Y$ are not fixed a
priori but adapted to the interaction of the two maps $P$ and $Q$ through
the intrisic Liouville measure. We note also that Theorem~\ref{thm:2} deals
with $L^2$ function spaces, while Theorem~\ref{thm:1} applies to general
direct integral decompositions of Hilbert spaces, as for instance to
splittings $H= \oplus_\lambda E_\lambda$ with non constant $\dim E_\lambda$
coming from spectral resolution of operators. See examples in
\S\ref{sec:applications-homogeneous}.

Another feature of Theorem~\ref{thm:1} is its invariance through general
isometries of $H$. Indeed, if an isometry $U$ acts on $H$, then two initial
POVM $P$ and $Q$ are conjugated to $P^U = U^{-1} P U$ and $Q^U = U^{-1}Q
U$. The Liouville measure is preserved: $\mu_{P^U Q^U} = \tau(P^UQ^U) =
\tau(PQ) =\mu_{PQ}$, as the measures associated to $\rho$ and $\rho^U$:
$\nu_P^\rho = \nu_{P^U}^{\rho^U}$ and $\nu_Q^\rho = \nu_{Q^U}^{\rho^U}$. In
comparison, \eqref{eq:4} is not invariant since the underlying isometry
$\mathcal{U}_{FL}$ there becomes $U^{-1}\mathcal{U}_{FL}U$, whose
$L^1\rightarrow L^\infty$ norm is a priori not controlled by
$\|\mathcal{U}_{FL}\|_{1,\infty}$ alone.

\section{Proof of Theorem~\ref{thm:1} and related comments}
\label{sec:proof-theorem}

Our proof of Theorem~\ref{thm:1} follows the lines of Frank and Lieb's
argument for Theorem~\ref{thm:2}. It relies on two classical lemmas in
quantum statistical mechanics, see e.g.~\cite{Carlen}
and~\cite[Thm.~8.5]{Simon}.

\begin{lemma}[Gibbs variational formula]
  \label{lemma:Gibbs}
  Let $A$ be a self-adjoint operator such that $e^{-A}$ is trace
  class. Then for any unit trace state $\rho$, it holds
  \begin{displaymath}
    \tau(A\rho) - S(\rho) \geq -\ln \tau(e^{-A})\,,
  \end{displaymath}
  with equality iff $\rho= e^{-A}/ \tau(e^{-A})$.
\end{lemma}

\begin{lemma}[Golden--Thompson inequality]
  \label{lemma:Golden-Thompson}
  Let $A$ and $B$ be self-adjoint operators with upper bounds. Then it
  holds that
  \begin{displaymath}
    \tau(e^{A+B}) \leq \tau(e^{A/2} e^{B} e^{A/2}) \,.
  \end{displaymath}
  Here, following \cite{Kato}, $A+B$ stands for the self-adjoint operator
  associated to the closed upper-bounded quadratic form $q_A + q_B$ on the
  closure of $D(A_{-}^{1/2})\cap D(B_{-}^{1/2})$.
\end{lemma}

We suppose that $\mu_{PQ}\leq \mu_P \otimes \mu_Q$ and first check the
absolute continuity of $\nu_P^\rho$ and $\nu_Q^\rho$ with respect to
$\mu_P$ and $\mu_Q$. If $\mu_P(\Omega)= 0$, then by $\sigma$-finiteness of
$\mu_Q$ one has $\mu_{PQ}(\Omega \times Y) = \tau(P(\Omega) ) =0$. Hence
$P(\Omega)= 0$, and $\nu_P^\rho(\Omega) = \tau(\rho P(\Omega))= 0$ as
claimed.

Given $k\geq 0$, let $ \varrho_{P,k} = \min (\varrho_P,k)$ and
$\varrho_{Q,k}= \min (\varrho_Q,k)$. We consider the following bounded
positive operators on $H$
\begin{displaymath}
  A_P = \int_X  \varrho_{P,k} d P
  \quad \mathrm{and} \quad
  A_Q = \int_Y  \varrho_{Q,k} d Q \,.
\end{displaymath}
Using a Hilbert basis $(e_i)$ of $H$ and monotone convergence, one sees
that
\begin{displaymath}
  \tau(A_P^{1/2}A_Q A_P^{1/2}) = \sum_i \langle A_Q A_P^{1/2}e_i ,
  A_P^{1/2} e_i \rangle = \int_Y  \varrho_{Q,k}
  \tau(A_P^{1/2} d Q A_P^{1/2}) \,, 
\end{displaymath}
where for any measurable set $\Omega' \subset Y$
\begin{align*}
  \tau(A_P^{1/2} Q(\Omega') A_P^{1/2}) &=
  \tau(Q^{1/2}(\Omega')A_P Q^{1/2}(\Omega')) \\
  & = \int_X \varrho_{P,k} \tau(Q^{1/2}(\Omega')d P Q^{1/2}(\Omega')) \\
  & = \int_X \varrho_{P,k} d \mu_{P,Q}( \cdot \times \Omega')\,.
\end{align*}
Hence we obtain
\begin{align}
  \tau(A_P^{1/2}A_Q A_P^{1/2}) & = \int_{X\times Y} \varrho_{P,k}
  \varrho_{Q,k} d \mu_{PQ} \leq \int_{X\times Y} \varrho_P \varrho_Q d
  \mu_{PQ}
  \nonumber \\
  & = \int_{X\times Y} \frac{d \nu_P^\rho}{d \mu_P} \frac{d \nu_Q^\rho}{d
    \mu_Q} d \mu_{PQ}
  \nonumber \\
  & \leq \int_{X\times Y} d\nu_P^\rho \otimes d\nu_Q^\rho
  \label{eq:5}\\
  & = \nu_P^\rho(X) \nu_Q^\rho (Y)= \tau(\rho)^2 = 1.  \nonumber
\end{align}
Then by Lemma~\ref{lemma:Golden-Thompson}, it holds that
\begin{equation}
  \label{eq:6}
  \tau(e^{\ln A_P + \ln A_Q})
  \leq \tau(A_P^{1/2} A_Q A_P^{1/2}) \leq 1\,,
\end{equation}
where
\begin{displaymath}
  A = \ln A_P +\ln A_Q 
\end{displaymath}
is the unique self-adjoint operator associated to the closed quadratic form
$q_{\ln A_P} + q_{\ln A_Q}$ on the closure $V$ of $D(A) =
D(\ln_{-}^{1/2}A_P) \cap D(\ln_{-}^{1/2}A_Q)$; see~\cite{Kato} and
\cite[VIII.6]{Reed-Simon}

We suppose now that $\rho$ writes $\rho = \sum_i p_i \Pi_{f_i}$ with
orthonomal functions $f_i$, and satisfies
\begin{displaymath}
  \int_X \ln_{-}\varrho_P d\nu_P^\rho < \infty \quad \mathrm{and}\quad
  \int_Y \ln_{-}\varrho_Q d\nu_Q^\rho < \infty \,.
\end{displaymath}
Note that otherwise \eqref{eq:1} is already satisfied. We show that $f_i
\in D(A)$ for any $i$, and evaluate $\langle A f_i,f_i\rangle$. We can use
a Jensen-type inequality for operator convex functions on POVM due to Choi,
see~\cite{Choi,Petz}. Namely, since $-\ln t$ is operator convex (see
\cite[Chap.V]{Bhatia}, it holds that
\begin{displaymath}
  -\ln A_P \leq -\int_X \ln \varrho_{P,k} d P\,,
\end{displaymath}
hence
\begin{align}
  \label{eq:7}
  - \langle \ln A_P f_i, f_i \rangle & \leq -\int_X \ln \varrho_{P,k} d
  \langle P f_i, f_i \rangle = -\int_X \ln \varrho_{P,k} d
  \nu_P^{\Pi_{f_i}} \\
  & \leq \int_X \ln_{-} \varrho_P d \nu_P^\rho < \infty \nonumber
\end{align}
by assumption on $\rho$. One gets a similar inequality for $ - \langle \ln
A_Q f_i, f_i \rangle$.  We can then apply Lemma~\ref{lemma:Gibbs} to $-A$
on $V$ with $\Pi_V \rho \Pi_V = \rho$. Summing \eqref{eq:7} and using
\eqref{eq:6} yields
\begin{equation}
  \label{eq:8}
  - \int_X \ln \varrho_{P,k} d \nu_P^\rho - \int_Y \ln \varrho_{Q,k} d \nu_Q^\rho
  \geq \tau(- A\rho) \geq S(\rho)\,,
\end{equation}
and gives Theorem~\ref{thm:1} by monotone convergence when $k \rightarrow
+\infty.$

\smallskip

\textbf{Equality case}. We discuss here the equality case in
Theorem~\ref{thm:1}. It holds iff \eqref{eq:5}, \eqref{eq:7}, \eqref{eq:8}
become sharp for $k\rightarrow \infty$. Equality in \eqref{eq:5} means that
the Liouville measure is a product, i.e.~ the POVM $P$ and $Q$ are
independent.  By \cite{Petz}, \eqref{eq:7} are all equalities iff $P$ and
$Q$ are \emph{projection} valued measures on the support of $\nu_P^\rho$
and $\nu_Q^\rho$. Then equality holds in \eqref{eq:8} iff it holds:
\begin{itemize}
\item in Golden-Thompson inequality, which is achieved when $P$ and $Q$
  commute,
\item and in Gibbs formula, meaning that $A_P A_Q \rightarrow \rho$,
  i.e. $\rho$ is a product state with respect to $P$, $Q$.
\end{itemize}

Altogether these conditions are very restrictive. Examples are given by a
space $H$ that splits into $H_1 \otimes H_2$ with given basis $(e_i)$,
$(f_j)$. There $P(i) = \Pi_{e_i} \otimes 1$ and $Q(j) = 1 \otimes
\Pi_{f_j}$, and product states $\rho= \rho_1 \otimes \rho_2$ give
equalities.

Note however that an approximate equality in \eqref{eq:1} may be achieved
in other cases. This happens for instance for the Fourier transform
discussed in \S\ref{sec:fourier-uncertainty}. Here the position-momentum
maps are independent projection valued measures, hence \eqref{eq:5} and
\eqref{eq:7} are equalities. Although $P$ and $Q$ don't commute here,
Theorem~\ref{thm:1} becomes sharp on states that spread at large scale on
the $(x,\xi)$ phase space like $\rho_t = e^{-t (\Delta + \|x\|^2)}/
\tau(e^{-t( \Delta + \|x\|^2)})$ when $t \searrow 0$; see \cite{Rumin11}.

\section{Applications on homogeneous spaces}
\label{sec:applications-homogeneous}

\subsection{Invariant operators and spatial/spectral independence}
\label{sec:invar-oper-indep}

We turn back to examples and describe a setting on homogeneous spaces
leading to independent spatial and spectral projections valued measures.

Let $H$ be a closed group in a locally compact group $G$ and let $X= G/H$
be the corresponding homogeneous space. We assume moreover that $X$ admits
a $G$-invariant Borel measure $\mu$. This is equivalent to the equality of
modular functions $\delta_H= \delta_G$ on $H$, and $\mu$ is unique up to
constants, see \cite{Helgason84,Nachbin}.

Let $V$ be a separable Hilbert space and consider $\mathcal{H}=
L^2(X,V,\mu)$. We say that an isometry $U$ of $\mathcal{H}$ is a gauge
transform if it acts fiberwise over $X$, i.e.~ $(Uf)(x) = U_x(f(x))$ for
some isometries $U_x:V\rightarrow V$.
\begin{defn}
  \label{defn:translation-invariance}
  An operator $A$ on $\mathcal{H}$ will be called \emph{translation
    gauge-invariant}, if it is conjugated by a gauge transform $U_0$ to an
  operator $A_0= U_0^{-1}A U_0$ whose translates $A_0^g= g^{-1}A_0g$ stay
  conjugated to $A_0$ up to gauge transforms, i.e.~ $A_0^g = U(g)^{-1} A_0
  U(g)$.
\end{defn}

This class contains the translation invariant operators on $X$, as
differential operators with constant coefficients when $G$ is a Lie
group. It contains also operators conjugated to translation invariant
operators by gauge transform; for instance $A = i\frac{d}{dx} + x =
e^{-ix^2/2} i\frac{d}{dx} e^{ix^2/2}$ is not translation-invariant, yet
translation gauge-invariant on $L^2(\R)$. Another classical example is
given by the Hamiltonian of a constant magnetic field $B \geq 0$ in $\R^2$
\begin{equation}
  \label{eq:9}
  H_B = \Bigl(-i\frac{\partial}{\partial x} + \frac{By}{2}\Bigr)^2 +
  \Bigl(-i\frac{\partial}{\partial y} - \frac{Bx}{2}\Bigr)^2 = (-i \nabla +
  A)^2\,. 
\end{equation}
A general gauge transform acts on $H_B$ as $e^{-if}H_B e^{if}= (-i\nabla
+A+df)^2$, and fixing $f(x,y) = \frac{B}{2}(y_0 x - x_0 y)$ actually
translates $H_B$ on $\R^2$ by $(x_0,y_0)$

Our interest for gauge invariant self-adjoint operators $A$ here comes from
the independence of the projection valued measure on $H$ coming from $X$;
namely $P(\Omega) = \chi_\Omega \times $, and the spectral resolution of
$A$ on $Y=\R$: $Q(I) = \Pi_A(I)$.

\begin{prop}
  \label{prop:independence}
  Let $A$ be a translation gauge-invariant self-adjoint operator on
  $L^2(X,V,\mu)$. Then there exists a measure $\mu_A$ on $\R$ such that the
  Liouville measure writes
  \begin{equation}
    \label{eq:10}
    \mu_{PQ} = \mu \otimes \mu_A \,.
  \end{equation}
\end{prop}

We shall call $\mu_A$ the \emph{spectral measure} of $A$ in the sequel (not
to be confused with the \emph{projection} valued measure $Q = \Pi_A$).
\begin{proof}
  We first observe that conjugating any $A$ by a gauge transform preserves
  $\mu_{PQ}$. Indeed if $A^U= U^{-1} A U$ then $\Pi_{A^U} (I) = U^{-1}
  \Pi_A (I) U$ and
  \begin{align*}
    \mu_{PQ^U} (\Omega \times I) & = \tau(\chi_\Omega \Pi_{A^U} (I) ) =
    \tau(
    U \chi_\Omega U^{-1} \Pi_A(I)) \\
    & = \tau (\chi_\Omega \Pi_A(I)) = \mu_{PQ}(\Omega\times I)\,.
  \end{align*}
  Then if $g^{-1} A g = A^{U(g)}$, one has
  \begin{align*}
    \mu_{PQ} (g.\Omega \times I) & = \tau(\chi_{g\Omega} \Pi_A (I)) =
    \tau(g \chi_\Omega g^{-1} \Pi_A (I) ) \\
    & = \tau (\chi_\Omega g^{-1} \Pi_A(I) g) = \tau (\chi_\Omega \Pi_{A^U}
    (I)) \\
    & = \mu_{PQ} (\Omega \times I)\,.
  \end{align*}
  Hence, given $I$, $\Omega \mapsto \mu_{PQ} (\Omega \times I)$ is an
  invariant measure on $X$ thus proportional to $\mu$ by uniqueness, i.e.~
  $\mu_{PQ}(\Omega \times I) = \mu(\Omega) \mu_A(I)$ as needed.
\end{proof}
More concretely, $\mu_A(I)$ expresses using the Schwarz kernel
$K_{\Pi_A(I)}$ of $\Pi_A(I)$. Indeed
\begin{align*}
  \mu_{PQ} (\Omega \times I) & = \tau_\mathcal{H}(\chi_\Omega \Pi_A(I) ) =
  \|\chi_\Omega \Pi_A(I)\|_{HS(\mathcal{H})}^2 \\
  &= \int_\Omega \int_X \|K_{\Pi_A(I)}(x,y)\|^2_{HS(V)} d \mu(y) d\mu(x),
\end{align*}
so that one has for $\mu$-almost every $x$
\begin{equation}
  \label{eq:11}
  \mu_A(I) = \int_X \|K_{\Pi_A(I)}(x,y)\|^2_{HS(V)} d \mu(y) \,. 
\end{equation}
This last formula, together with Plancherel formula on $X$, helps in
computing examples.

\subsection{Compact case}
\label{sec:compact-case}

We first consider the case of homogeneous spaces $X=G/H$ for compact groups
$G$. Then a scalar invariant self-adjoint operator $A$ on $X$ induces a
spectral splitting
\begin{equation}
  \label{eq:12}
  L^2(X)= \bigoplus_{\Sp A} E_\lambda\,.
\end{equation}
The spectral measure $\mu_A$ is supported on $\Sp A$ and
\begin{align*}
  \mu_{PQ}(X \times \Pi_\lambda) & = \tau(\Pi_\lambda) = \dim
  E_\lambda  \\
  & = \mu(X) \mu_A(\lambda)
\end{align*}
hence $\mu_A(\lambda) = \dim E_\lambda / \mu(X)$. Therefore \eqref{eq:1}
writes here
\begin{equation}
  \label{eq:13}
  -\int_X \varrho(x) \ln \varrho(x) d\mu(x) - \sum_{\Sp A} \ln \Bigl(
  \frac{\tau(\rho \Pi_\lambda)}{\dim E_\lambda}\Bigr) \tau(\rho \Pi_\lambda) \geq
  S(\rho)\,, 
\end{equation}
if $\mu$ is normalized such that $\mu(X) = 1$.

A first remark here is that the spectral entropy sum, say $S_A (\rho)$,
actually does not depend on the values in $\Sp A$. It stays the same for
any other operator $\varphi(A)$ with $\varphi$ injective on $\Sp A$. This
sum only depends on the repartition on the state in the splitting
\eqref{eq:12}.

This stays true in general for the change of self-adjoint operator $A$ into
$\varphi(A)$ on any space $X$. Indeed one has $\Pi_{\varphi(A)}(I) = \Pi_A
(\varphi^{-1}(I))$ and the induced projection valued measure becomes
$Q_{\varphi(A)} = Q_A \circ \varphi^{-1}$. One sees easily (or by
Proposition~\ref{prop:refinement} below) that this amounts in a change of
variables $\lambda \mapsto \varphi^{-1}(\lambda)$ in the spectral entropy
integrals over $Y= \R$ in \eqref{eq:1}, and thus $S_A(\rho)=
S_{\varphi(A)}(\rho)$. Hence this notion does not depend on the actual
values of energy levels, not even their order, but deals with the
repartition of the state in the direct integral splitting of $H
=\int^\oplus dE_\lambda$.

We also observe that \eqref{eq:13} is sharp on purely spectral states such
that $\rho=\varphi(A)$ with $\tau(\rho)= 1$. Indeed these translation
invariant states have constant density $\varrho(x) = 1$ since $\mu(X) = 1 =
\tau(\rho)$. Thus the spatial entropy term vanishes in \eqref{eq:13}, while
$$
S_A(\rho) = S(\rho) = -\sum_{\Sp A} (\ln \varphi(\lambda)) \varphi(\lambda)
\dim E_\lambda \,.
$$

Finally, we note that $S_Q(\rho)$ decreases with the refinement of the
splitting of $L^2(X)$ and is ultimately bounded from below by the finest
possible one coming from the decomposition of the unitary representation
$\pi$ of $G$ in $L^2(X)$ into finite dimensional irreducible
representations:
\begin{equation}
  \label{eq:14}
  L^2(X) = \bigoplus_{\sigma \subset \pi} H_\sigma \,.
\end{equation}
This comes from the splitting of each invariant space $E_\lambda$ into such
sums of $H_\sigma$ (explaining the discreteness of $\Sp A$), and the
following general monotony property.

\begin{prop}
  \label{prop:refinement}
  Let $P$ be a POVM from $(X,\mu)$ to $H$. Suppose that $P'$ is a POVM on
  $(X',\mu')$ induced from $P$ by a measurable map $\varphi: X \rightarrow
  X'$, i.e.~ $P'= P(\varphi^{-1})$ and $\mu'= \mu (\varphi^{-1})$. Suppose
  moreover that $X'$ is $\sigma$-finite. Then
  \begin{equation}
    \label{eq:15}
    S_{P'}(\rho) = - \int_{X'} \ln \bigl( \frac{d\nu_{P'}^\rho}{d \mu'} \bigr)
    d\nu_{P'}^\rho \geq S_P(\rho) = - \int_X \ln \bigl( \frac{d\nu_P^\rho}{ d\mu}
    \bigr) d \nu_P^\rho \,,
  \end{equation}
  provided these integrals are defined.
\end{prop}

Indeed, Jensen inequality applied to the conditional expectation $E\bigl(
\frac{d \nu_P^\rho}{d\mu} | \varphi \bigr) = \frac{d \nu_{P'}^\rho}{d
  \mu'}$ gives
\begin{displaymath}
  E\Bigl( - \frac{d \nu_P^\rho}{d\mu} \ln \bigl(\frac{d \nu_P^\rho}{d\mu}
  \bigr)| \varphi \Bigr) \leq - \frac{d \nu_{P'}^\rho}{d\mu'} \ln
  \bigl(\frac{d \nu_{P'}^\rho}{d\mu'} \bigr)\,,
\end{displaymath}
yielding \eqref{eq:15} by integration.
  

\smallskip

Eventually we illustrate this discussion on the sphere $S^{n-1} =
\mathrm{SO}(n)/ \mathrm{SO}(n-1)$ and the Laplacian $\Delta_S$. Here
$\Sp(\Delta_S)= \{d(d+n-2) \mid d\geq 0\}$ and $E_d(\Delta)$ consists in
the harmonic polynomials of degree $d$, with $\dim E_d = \binom{d+n-1}{n-1}
- \binom{d+n-3}{n-1}$, see e.g.~ \cite{Helgason84,Strichartz}. Moreover the
representation of $\mathrm{SO}(n)$ on $E_d$ is irreducible, so that the
spectral and Weyl decompositions \eqref{eq:12} and \eqref{eq:14} coincide
here. That means that the measurement of $\rho$ by $\Delta_S$ (its energy
distribution), actually provides the best (lowest) spectral entropic term
(among invariant operators) in \eqref{eq:13}.

\subsection{Non-compact examples}
\label{sec:non-compact-examples}

We compute the spectral measure of some operators in non-compact
situations. Let $A$ be a self-adjoint differential operator with constant
coefficients on $X=\R^n$ and $\sigma_A$ its polynomial symbol. The Fourier
transform of the spectral projection $\Pi_A(I)$ is the multiplication by
$\chi_{\sigma_A^{-1}(I)}$ on $\R^n_\xi$. Then by \eqref{eq:11} and
Plancherel formula the spectral measure of $A$ is
\begin{equation}
  \label{eq:16}
  \mu_A(I) = \|\chi_{\sigma^{-1}(I)}\|_{L_\xi^2}^2 = \mathrm{vol}
  (\sigma_A^{-1}(I)) \,. 
\end{equation}
Using the coarea formula, its density with respect to Lebesgue measure on
$\R$ is
\begin{displaymath}
  \frac{d \mu_A}{d\lambda} = \int_{\sigma_A^{-1}(\lambda)} \frac{d
    H_\sigma}{|\nabla \sigma|}
\end{displaymath}
for the hypersurface measure $H_\sigma$ on
$\sigma^{-1}(\lambda)$. Similarly the measure of a state $\rho = \sum_i p_i
\Pi_{f_i}$ in $L^2(\R^n)$ relatively to the spectral resolution of $A$
reads
\begin{displaymath}
  \nu^\rho_A (I) = \tau( \rho \Pi_A(I)) = \int_{\sigma^{-1}_A(I)}
  \widehat\varrho(\xi) d\xi   
\end{displaymath}
where $ \widehat\varrho(\xi) = \sum_i p_i |\widehat f_i(\xi)|^2$ as in
\S\ref{sec:fourier-uncertainty}. Hence
\begin{displaymath}
  \frac{d\nu_A^\rho}{d\mu_A}(\lambda) = \int_{\sigma_A^{-1}(\lambda)}
  \widehat\varrho(\xi) \frac{d
    H_\sigma}{|\nabla \sigma|} / \int_{\sigma_A^{-1}(\lambda)} \frac{d
    H_\sigma}{|\nabla \sigma|}\,,
\end{displaymath}
expliciting the spectral entropy term in \eqref{eq:1}:
\begin{equation}
  \label{eq:17}
  S_A(\rho) = -\int_{\R} \ln \bigl(\frac{d\nu^\rho_A}{d \mu_A}\bigr) d
  \nu_A^\rho \,.
\end{equation}

For instance in the case of the Laplacian $\Delta$. Then
$\sigma_\Delta(\xi) = \|2 \pi\xi\|^2$ with our convention in
\S\ref{sec:fourier-uncertainty}, and the spectral measure is
$\nu_\Delta([0,\lambda]) = \mathrm{vol}(B(\sqrt \lambda / 2\pi))$. Then
$H_\sigma$ is the usual measure $m$ on the spheres $S_R$, and one finds
that
\begin{displaymath}
  S_\Delta(\rho) =- \int_0^\infty \ln \Bigr(\int_{S_R} \widehat
  \varrho(\xi) \frac{d m(\xi)}{m(S_R)} \Bigl) \Bigl(\int_{S_R} \widehat
  \varrho(\xi) d m(\xi) \Bigr) dR \
\end{displaymath}

We note that for any invariant $A$, the map $\Pi_A$ is actually induced (up
to Fourier transform) from the momentum map $Q(\Omega) = \chi_\Omega \times
$ on $ \R^n_\xi$ by the symbol $\sigma_A :\R^n_\xi \rightarrow \R$; i.e.~
$\Pi_A = Q(\sigma_A^{-1})$ and $\mu_A = \mu_Q (\sigma_A^{-1})$ by
\eqref{eq:16}. Hence Proposition~\ref{prop:refinement} applies and gives
the lower bound
\begin{displaymath}
  S_A(\rho) \geq S_Q (\rho) = - \int_{\R^n} \widehat \varrho(\xi) \ln
  \widehat\varrho(\xi)  d\xi \,,
\end{displaymath}
arising in the Fourier uncertainty principle \eqref{eq:3}. Yet, contrarily
to the Laplacian on $S^{n-1}$, one has in general $S_A(\rho) > S_Q(\rho)$,
as due to Jensen inequality. Indeed the knowledge of $S_A(\rho)$ requires
less information on $\rho$ than $S_Q(\rho)$, since $\nu^\rho_A$ only
depends on the mean values of $\widehat\varrho$ on the level sets of
$\sigma_A$.

\smallskip

We mention that the spectral measure $\mu_A$ can (in principle) be computed
on other non-compact homogeneous spaces $X=G/H$, as long an explicit
Plancherel formula is available. This is indeed the case for symmetric
spaces $G/K$ with $G$ connected semi-simple Lie group with finite center
and $K$ its maximal compact subgroup; see e.g.~\cite{Helgason84}.
For instance, computations of the spectral measure of the Laplacian on
symmetric spaces may be found in \cite{Strichartz}.

\smallskip

We close this series with the example \eqref{eq:9} of the Hamiltonian
$H_B$ in $L^2(\R^2)$. It is translation gauge-invariant and
Proposition~\ref{sec:invar-oper-indep} applies. Its spectrum is discrete
and consists in the (Landau levels) $\lambda_n= (2n+1)B$, $n\in \N$, with
each eigenspace $E_n$ of constant density $\mu_{H_B}(\lambda_n) =
\frac{B}{2\pi}$; see e.g.~\cite{ELV}. Hence \eqref{eq:1} reads here
\begin{equation}
  \label{eq:18}
  -\int_{\R^2} \varrho(x) \ln \varrho(x) dx - \sum_{n\geq 0} 
  \tau(\rho \Pi_n) \ln \tau(\rho \Pi_n) \geq S(\rho) - \ln (B/2\pi)\,.
\end{equation}
One sees that the uncertainty constrain relaxes for a given state when $B$
increases. Indeed the state may concentrate on fewer Landau levels, whose
density increases with $B$. Note also that the spectral entropy term in the
left side is always positive here and vanishes iff $\rho$ is contained in a
single level, in which case $-\int_{\R^2} \varrho(x) \ln \varrho(x) dx \geq
S(\rho) - \ln (B/2\pi)$.

\subsection{Log-Sobolev inequalities}
\label{sec:log-sobol-ineq-1}

As emphasized above, the spectral entropy associated to an operator does
not depend on the actual values of energy levels, but only on the induced
spectral splitting. It turns out however that one can bound it using a
single estimation of the mean energy of the state $\mathcal{E}_A(\rho) =
\tau(A \rho)$; yielding a log-Sobolev (entropy-energy) inequality.

Consider again a self-adjoint translation gauge-invariant operator $A$ on
an homogeneous space $X=G/H$. For $t\geq 0$, let
\begin{equation}
  \label{eq:19}
  L_A(t) = \int_\R e^{-t\lambda} d\mu_A (\lambda)
\end{equation}
denotes the Laplace transform of the spectral measure $\mu_A$. This is
actually also the constant ratio $\tau(e^{-tA} \chi_\Omega)/ \mu(\Omega)$,
or heat decay of $e^{-tA}$. The following Gibbs inequality holds.
\begin{prop}
  \label{prop:entropy-laplace}
  Suppose that the state $\rho$ satisfies $\mathcal{E}_A^+(\rho) =\tau(\max
  (A,0) \rho) < \infty$ and that $L(t) < \infty$. Then
  \begin{displaymath}
    \int_\R \ln^- \bigl( \frac{d\nu_A^\rho}{d\mu_A} \bigr) d\mu_A(\rho) < \infty 
  \end{displaymath}
  and it holds that
  \begin{equation}
    \label{eq:20}
    S_A(\rho) \leq t \mathcal{E}_A(\rho) + \ln L_A(t) \,.
  \end{equation}
\end{prop}

\begin{proof}
  From $\mathcal{E}_A(\rho) = \tau (\rho A) = \int_\R \lambda \tau(\rho d
  \Pi_A) = \int_\R \lambda d\nu_A^\rho$, one has
  \begin{align*}
    -\int_\R \bigl( \frac{d\nu_A^\rho}{e^{-t\lambda} d\mu_A} \bigr) \ln
    \bigl( \frac{d\nu_A^\rho}{e^{-t\lambda} d\mu_A} \bigr) e^{-t\lambda} d
    \mu_A & = S_A(\rho) - t \int_\R \lambda d\nu_A^\rho  \\
    &  =S_A(\rho) -t \mathcal{E}_A(\rho) \\
    & \leq -\Bigl( \int_\R d\nu_A^\rho\Bigr) \ln \Bigl( \int_\R \frac{d
      \nu_A^\rho}{L(t)} \Bigr) = \ln L(t) \,,
  \end{align*}
  by Jensen inequality and $\int_\R d\nu_A^\rho\ = \nu_A^\rho(\R) =
  \tau(\rho) =1$. Note that $ \ln^-(\frac{d\nu_A^\rho}{d\mu_A}) -t
  \lambda^+ \leq (-\ln)^+ (\frac{d\nu_A^\rho}{e^{-t\lambda} d\mu_A})$ whose
  corresponding integral is finite by Jensen if $L_A(t) < \infty$. This
  ensures that $S_A(\rho)$ exists.
\end{proof}

\begin{rem}
  \label{rem:Gibbs}
  Note moreover that equality holds iff $d\nu_A^\rho = C e^{-t\lambda} d
  \mu_A = e^{-t\lambda} d \mu_A / L_A(t)$, i.e. when the spectral
  distribution of the state is a Gibbs measure, i.e.~with exponential law
  with respect to energy. Such states always exist in the homogeneous case:
  take for instance
  \begin{displaymath}
    \rho_t =  e^{-\frac{t}{2}A} \chi_\Omega e^{-\frac{t}{2}A} / L_A(t) \mu(\Omega) \,.
  \end{displaymath}
\end{rem}

A corollary is the following log-Sobolev inequality, in a parametric form.
\begin{cor}
  \label{cor:log-sob}
  Under the previous assumptions, it holds that
  \begin{equation}
    \label{eq:21}
    -\int_X \varrho(x) \ln \varrho(x) d\mu(x) + t \mathcal{E}_A(\rho) + \ln
    L_A(t)  \geq S(\rho) \,.
  \end{equation}
\end{cor}

When $\rho$ is a pure state $\Pi_f$ and $A$ is the Laplacian (or $e^{-tA}$
is Markovian) this is a well known inequality, see \cite{Davies}. It
extends here for mixed states on a large class of operators: translation
gauge-invariant ones on homogeneous manifolds.

One can also optimize \eqref{eq:21} in $t$. Notice that $\ln L_A(t)$ is a
convex function, since $L_A(t)$ is the integral of log-convex functions. We
consider then the Legendre transform (Young conjugate) of $-\ln L_A(t)$,
namely
\begin{displaymath}
  (\ln L_A)^*(\lambda) = \inf_{t\geq 0} (t \lambda + \ln L_A(t)) \,.
\end{displaymath}
Then \eqref{eq:21} reads
\begin{equation}
  \label{eq:22}
  -\int_X \varrho(x) \ln \varrho(x) d\mu(x) + (\ln
  L_A)^*(\mathcal{E}_A(\rho)) \geq
  S(\rho) \,.
\end{equation}
Note that by Remark~\ref{rem:Gibbs}, \eqref{eq:22} is also equivalent to
the previous purely entropic inequality
\begin{displaymath}
  -\int_X \varrho(x) \ln \varrho(x) d\mu(x) + S_A(\rho) \geq S(\rho)
\end{displaymath}
on states at statistical equilibrium with respect to the energy $A$, i.e.~
with spectral measure $d\nu_A^\rho = C e^{-t\lambda} d \mu_A$ for some $t$,
but is weaker on more general states. Moreover, following the discussion in
\S\ref{sec:compact-case}, the inequality \eqref{eq:22} is sharp (on mixed
states) on homogeneous spaces $X= G/H$ with $G$ compact, if $e^{-tA}$ is
trace class for some $t$.

\smallskip

\textbf{Examples.}  In the case of the Laplacian on $\R^n$, one has
$L_\Delta(t) = (4\pi t)^{-n/2}$; a well known heat-decay (that also follows
from \eqref{eq:19} with $\nu_\Delta([0,\lambda])$ given in
\S\ref{sec:non-compact-examples}). Then \eqref{eq:22} writes
\begin{equation}
  \label{eq:23}
  -\int_X \varrho(x) \ln \varrho(x) d\mu(x) + \frac{n}{2} \ln \Bigl(
  \frac{e \mathcal{E}_\Delta(\rho)}{2\pi n} \Bigr) \geq S(\rho) \,.
\end{equation}
This inequality is due to Dolbeault-Felmer-Loss-Paturel in \cite{DFLP} and
asymptotically sharp on normalized harmonic oscillator heat $\rho_t = e^{-t
  (\Delta + \|x\|^2)}/ \tau(e^{-t (\Delta + \|x\|^2)})$ for $t \searrow 0$.

\smallskip

We can compare this on $\R^2$ with the Hamiltonian $H_B$ in
\eqref{eq:9}. From discussion in \S\ref{sec:non-compact-examples}, one has
here
\begin{displaymath}
  L_B(t) = \sum_{n\geq 0} \frac{B}{2 \pi} e^{-(2n +1)Bt} = \frac{B}{4\pi
    \sinh(Bt)} \,.
\end{displaymath}
Straightforward calculations then yields
\begin{equation}
  \label{eq:24}
  (\ln L_B)^*(\mathcal{E}_{H_B}(\rho)) = \ln \Bigl( 
  \frac{\mathcal{E}_{H_B}(\rho) + B}{4\pi} \Bigr) + n(\rho) \ln \bigl(1 +
  \frac{1}{n(\rho)} \bigr) \,,
\end{equation}
where
\begin{displaymath}
  n(\rho) = \sum_n n \tau(\Pi_n \rho)
\end{displaymath}
is the mean Landau level of $\rho$ (with respect to the ground state). Note
that
\begin{displaymath}
  \mathcal{E}_{H_B}(\rho) = (2n(\rho) + 1) B \,.
\end{displaymath}
When $B\rightarrow 0$, one has $H_B \rightarrow \Delta $ and $n(\rho)
\rightarrow + \infty$ for a given state. Then
\begin{displaymath}
  (\ln L_B)^*(\mathcal{E}_{H_B}(\rho)) \rightarrow \ln \bigl( \frac{e
    \mathcal{E}_\Delta(\rho)}{4\pi} \bigr)
\end{displaymath}
and one recovers \eqref{eq:23} on $\R^2$ from \eqref{eq:22}.


\smallskip

\textbf{Further comments.}  The log-Sobolev inequality \eqref{eq:22}
obtained here may be compared to an other similar result proved in
\cite{Rumin11}. There Corollary 1.6 states that
\begin{equation}
  \label{eq:25}
  -\int_X \varrho(x) \ln \varrho(x) d\mu(x) + (\ln
  F_A)^c(\mathcal{E}_A(\rho)) \geq -3 - \ln \|\rho\|_{L^2 \rightarrow L^2} 
\end{equation}
where
\begin{displaymath}
  F_A (\lambda)= \sup_\Omega \Bigl(\frac{\tau(\chi_\Omega \Pi_A(]-\infty,
    \lambda[)}{\mu(\Omega)} \Bigr)
\end{displaymath}
and $(\ln F_A)^c$ denotes the concave hull of $\ln F_A$. This statement
holds on general $\sigma$-finite spaces without invariance and homogeneity
assumptions.

In an homogenous situation, one has $F_A(\lambda) = \mu_A(]-\infty,
\lambda[)$. It turns out that
\begin{equation}
  \label{eq:26}
  (\ln L_A)^* \geq (\ln F_A)^c
\end{equation}
so that the left side of \eqref{eq:22} is larger than the one in
\eqref{eq:25}. Indeed,
\begin{align*}
  L_A(t) & = \int_\R e^{-tu} d\mu_A (u) \\
  & \geq \int_{-\infty}^\lambda e^{-t\lambda} d\mu_A(u) = e^{-t\lambda}
  F_A(\lambda)
\end{align*}
and thus $\ln F_A(\lambda) \leq t \lambda + \ln L_A(t)$, giving
\eqref{eq:26} by concavity of $(\ln L_A)^*$. On the other hand the right
side of \eqref{eq:25} is smaller than in \eqref{eq:22} since on unit trace
states
\begin{displaymath}
  -\ln \|\rho\|_{L^2\rightarrow L^2} \leq S(\rho) = -\tau(\rho \ln \rho) \,,
\end{displaymath}
with equality on uniformly distributed states, i.e.~normalized projections
on finite dimensional subspaces of $H$.

In conclusion the two log-Sobolev inequalities we discuss here are not
equivalent, even on homogeneous spaces; see also \cite[\S4.2]{Rumin11} for
a more precise comparison in the case of the Laplacian on $\R^n$. Yet, we
have seen that the version developped here is sharp in some classical
cases, including translation gauge-invariant operators on compact
homogeneous spaces. Moreover it comes from the stronger entropic
uncertainty principle stated in Theorem~\ref{thm:1}.

\smallskip

\textbf{Acknowledgments}. The author is grateful to Rupert L. Frank and
Elliott H. Lieb for showing him their earlier proof of the Fourier
uncertainty inequality \eqref{eq:3}.

\bibliographystyle{abbrv}



\begin{thebibliography}{10}

\bibitem{Bhatia}
R.~Bhatia.
\newblock {\em Matrix analysis}, volume 169 of {\em Graduate Texts in
  Mathematics}.
\newblock Springer-Verlag, New York, 1997.

\bibitem{Carlen}
E.~Carlen.
\newblock Trace inequalities and quantum entropy: an introductory course.
\newblock In {\em Entropy and the quantum}, volume 529 of {\em Contemp. Math.},
  pages 73--140. Amer. Math. Soc., Providence, RI, 2010.

\bibitem{Choi}
M.~D. Choi.
\newblock A {S}chwarz inequality for positive linear maps on {$C^{\ast} \
  $}-algebras.
\newblock {\em Illinois J. Math.}, 18:565--574, 1974.

\bibitem{Davies}
E.~B. Davies.
\newblock {\em Heat kernels and spectral theory}, volume~92 of {\em Cambridge
  Tracts in Mathematics}.
\newblock Cambridge University Press, Cambridge, 1989.

\bibitem{DFLP}
J.~Dolbeault, P.~Felmer, M.~Loss, and E.~Paturel.
\newblock Lieb-{T}hirring type inequalities and {G}agliardo-{N}irenberg
  inequalities for systems.
\newblock {\em J. Funct. Anal.}, 238(1):193--220, 2006.

\bibitem{ELV}
L.~Erd{\H{o}}s, M.~Loss, and V.~Vougalter.
\newblock Diamagnetic behavior of sums of {D}irichlet eigenvalues.
\newblock {\em Ann. Inst. Fourier (Grenoble)}, 50(3):891--907, 2000.

\bibitem{Frank-Lieb}
R.~L. Frank and E.~H. Lieb.
\newblock Entropy and the uncertainty principle.
\newblock Preprint arXiv \url{http://fr.arxiv.org/abs/1109.1209}.

\bibitem{Helgason84}
S.~Helgason.
\newblock {\em Groups and geometric analysis}, volume 113 of {\em Pure and
  Applied Mathematics}.
\newblock Academic Press Inc., Orlando, FL, 1984.
\newblock Integral geometry, invariant differential operators, and spherical
  functions.

\bibitem{Kato}
T.~Kato.
\newblock Trotter's product formula for an arbitrary pair of self-adjoint
  contraction semigroups.
\newblock In {\em Topics in functional analysis (essays dedicated to {M}. {G}.
  {K}re\u\i n on the occasion of his 70th birthday)}, volume~3 of {\em Adv. in
  Math. Suppl. Stud.}, pages 185--195. Academic Press, New York, 1978.

\bibitem{Nachbin}
L.~Nachbin.
\newblock {\em The {H}aar integral}.
\newblock Robert E. Krieger Publishing Co., Huntington, N.Y., 1976.
\newblock Translated from the Portuguese by Lulu Bechtolsheim, Reprint of the
  1965 edition.

\bibitem{Nielsen-Chuang}
M.~A. Nielsen and I.~L. Chuang.
\newblock {\em Quantum computation and quantum information}.
\newblock Cambridge University Press, Cambridge, 2000.

\bibitem{Petz}
D.~Petz.
\newblock On the equality in {J}ensen's inequality for operator convex
  functions.
\newblock {\em Integral Equations Operator Theory}, 9(5):744--747, 1986.

\bibitem{Reed-Simon}
M.~Reed and B.~Simon.
\newblock {\em Methods of modern mathematical physics. {I}}.
\newblock Academic Press Inc. [Harcourt Brace Jovanovich Publishers], New York,
  second edition, 1980.
\newblock Functional analysis.

\bibitem{Rumin11}
M.~Rumin.
\newblock Balanced distribution--energy inequalities and related entropy
  bounds.
\newblock Preprint 2010. To appear in Duke Math Journal.

\bibitem{Simon}
B.~Simon.
\newblock {\em Trace ideals and their applications}, volume 120 of {\em
  Mathematical Surveys and Monographs}.
\newblock American Mathematical Society, Providence, RI, second edition, 2005.

\bibitem{Strichartz}
R.~S. Strichartz.
\newblock Estimates for sums of eigenvalues for domains in homogeneous spaces.
\newblock {\em J. Funct. Anal.}, 137(1):152--190, 1996.

\bibitem{Wikipedia-MUB}
Wikipedia.
\newblock Mutually unbiased bases.
\newblock \url{http://en.wikipedia.org/wiki/Mutually_unbiased_bases}.

\bibitem{Wikipedia-POVM}
Wikipedia.
\newblock P{O}{V}{M}.
\newblock \url{http://en.wikipedia.org/wiki/Positive_operator-valued_measure}.

\end{thebibliography}

\def\cprime{$'$} \def\dbar{\leavevmode\hbox to 0pt{\hskip.2ex
  \accent"16\hss}d}

-----------------------------------------------------

\end{document}